\documentclass{article}
\usepackage{amsmath}
\usepackage{amsfonts}
\usepackage{amsthm}
\usepackage{geometry}
\title{Foundations of abstract probability theory}
\author{Yurii Yurchenko\thanks{Odessa Polytechnic National University, Institute of Computer Systems, Department of Applied Mathematics and Information Technology, Shevchenko av. 1, Odesa 65044, Ukraine}}
\date{}
\begin{document}
\maketitle
\begin{abstract}
Using the ideas of abstract algebra, we introduce the basic concepts of abstract probability theory that generalize the Kolmogorov's probability theory, possibility theory and other theories that deal with uncertainty. Based on abstract addition and multiplication, we define an abstract measure and abstract Lebesgue integral. System of Kolmogorov's axioms is criticized, after which we introduce an abstract probability measure and abstract conditional probability, show that they have recognizable probability properties. In addition, we define an abstract expected value operator as the abstract Lebesgue integral and prove its properties. 
\end{abstract}
\newtheorem{ax}{Axiom}
\newtheorem{axr}{A}
\newtheorem{axkol}{Kolmogorov's axiom}
\newtheorem{cor}{Corollary}[section]
\newtheorem{prop}{Proposition}[section]
\newtheorem{defi}{Definition}[section]
\newtheorem{property}{Property}[section]
\newtheorem{lemma}{Lemma}[section]
\newtheorem{theorem}{Theorem}[section]
\section{Introduction}
The fundamentals of classical probability theory were formulated by Kolmogorov \cite{kol}, he axiomatically defined the theory of probability using the concepts of measure theory \cite{halmos}. The axioms are also supplemented by two definitions of conditional probability and independent events. In addition, for an real-valued random variable the expected value of the random variable is defined as the Lebesgue integral \cite{burkill}.

As an alternative to probability theory, was introduced the possibility theory that also deals with uncertainty. The possibility measure, as the membership function, was originally proposed by Zadeh \cite{zadeh} to describe his theory of fuzzy set. Later, this theory was studied by Dubois and Prade \cite{dub1,dub2}, who formed it in the form in which we know it now.

The difference between these theories is that in the probability theory the probability measure is $+$-additive, while in the possibility theory the possibility measure is $\vee$-additive, where $\vee$ is a logic sum, i.e. maximum of two values.

In addition, there are other generalizations of the concept of porobability.  Feynman \cite{feynman} and Burgin \cite{burgin} considered in their works the concept of a negative probability. Other examples of such a generalization are probability theories where the probability is hyperbolic number \cite{hyperbolic}.

Summng up all of the above, it becomes quite obvious that the concept of probability can be defined for different algebraic structures.

The concept of algebraic structure is one of the basic concepts in abstract algebra. According to Grillet \cite{grillet}, an algebraic structure can be defined as a set equipped with operations called addition and multiplication. Examples of algebraic structures are monoids, groups, rings and fields.

Thus, the classical probability is defined for an algebraic structure $(\mathbb{R},+,\times)$, which is a field, while the possibility is defined for an algebraic structure $(\mathbb{R},\vee,\wedge)$, which is a semiring with commutative multiplication. Moreover, note that, as we described earlier, probability not necessarily take real values, for exanple, hyperbolic probability is defined for an algebraic structure $(\mathbb{D},+,\times)$.

The aim of our work is to create an axiomatics of a theory that would generalize all the above examples of theories that deal with uncertainty.
\section{Abstract probability as abstract measure}
The concept of measure \cite{halmos} is fundamental in measure theory and is widely used in various areas of mathematics, for example, the axiomatics of Kolmogorov's probability theory is largely based on concepts presented in measure theory. It is clear that the algebraic structure $(\mathbb{R},+,\times)$ plays a key role in the definition of the classical measure.

This concept  was also attempted to be generalized to other algebraic structures. For example, Sugeno introduced the fuzzy measure \cite{sugeno}, which can be considered as a type of measure defined on algebraic structure $([0,1]\subset\mathbb{R},\vee,\wedge)$. In addition, many authors expressed their own vision of the fuzzy measure:  pseudo-additive measure \cite{sugeno1}  defined on  $([0,+\infty]\subset\mathbb{R},\oplus,\otimes)$, where $\oplus$ is pseudo-addition and $\otimes$ is pseudo-multiplication, and also \cite{guo},  \cite{ichi}, \cite{traore}, \cite{zhang}, \cite{zhang1}. All the measures described above are united by fact that they take values of non-negative real numbers and additive with respect to specific operations instead of the usual addition. 

Summarizing the above, we can distiguish the properties of a generalized algebraic structure.
\begin{defi}
\label{algstr}
Algebraic structure $(\mathbb{S},\oplus,\otimes,0_\mathbb{S},1_\mathbb{S},\succeq,d)$ is a metric partially ordered pseudo-semiring if it satisfies the following properties:
\begin{enumerate}
\item $(\mathbb{S},\succeq)$ is a non-strict partially ordered set;
\item Set $\mathbb{S}$ equipped with a binary operations $\oplus:\mathbb{S}\times \mathbb{S}\to \mathbb{S}$ and $\otimes:\mathbb{S}\times \mathbb{S}\to \mathbb{S}$ that satisfy the following properties:
\begin{itemize}
\item Commutativity of addition: $\forall a,b\in\mathbb{S}\;(a\oplus b=b\oplus a)$;
\item Commutativity of multiplication: $\forall a,b\in\mathbb{S}\;(a\otimes b=b\otimes a)$;
\item Associativity of addition: $\forall a,b,c\in\mathbb{S}\;(a\oplus(b\oplus c)=(a\oplus b)\oplus c)$;
\item Associativity of multiplication: $\forall a,b,c\in\mathbb{S}\;(a\otimes(b\otimes c)=(a\otimes b)\otimes c)$;
\item Additive identity $0_\mathbb{S}$: $\exists 0_\mathbb{S}\in\mathbb{S}:\forall a\in\mathbb{S}\; (a\oplus 0_\mathbb{S}=0_\mathbb{S}\oplus a=a)$;
\item Multiplicative identity $1_\mathbb{S}$: $\exists 1_\mathbb{S}\in\mathbb{S}:\forall a\in\mathbb{S}\; (a\otimes 1_\mathbb{S}=1_\mathbb{S}\otimes a=a)$;
\item Monotonicity of addition: $\forall a\in\mathbb{S}\; (b\succeq c\Rightarrow a\oplus b\succeq a\oplus c)$;
\item Monotonicity of multiplication: $\forall a\succeq 0_\mathbb{S}\; (b\succeq c\Rightarrow a\otimes b\succeq a\otimes c)$;
\item Distributivity of multiplication over addition: $\forall a,b,c\in\mathbb{S}\;(a\otimes(b\oplus c)=(a\otimes b)\oplus(a\otimes c))$;
\end{itemize}
\item $(\mathbb{S},d)$ is a metric space that satisfy the following properties:
\begin{itemize}
\item If $\{a_n\}$ and $\{b_n\}$ are sequences that $\forall n(a_n\preceq b_n)$, $d(a_n,a)\to 0$ and $d(b_n,b)\to 0$ as $n\to\infty$, then $a\preceq b$;
\item If $\{a_n\}$ and $\{b_n\}$ are sequences that $d(a_n,a)\to 0$ and $d(b_n,b)\to 0$ as $n\to\infty$, then $d(a_n\oplus b_n,a\oplus b)\to 0$ as $n\to\infty$;
\item If $\{a_n\}$ and $\{b_n\}$ are sequences that $d(a_n,a)\to 0$ and $d(b_n,b)\to 0$ as $n\to\infty$, then $d(a_n\otimes b_n,a\otimes b)\to 0$ as $n\to\infty$;
\end{itemize}
\end{enumerate}
\end{defi}
Here and below, $(\mathbb{S},\oplus,\otimes,0_\mathbb{S},1_\mathbb{S},\succeq,d)$ is always the metric partially ordered pseudo-semiring.
\begin{defi}
\label{meas}
Let $\left(\Omega, \mathcal{F}\right)$ is a measurable space. Abstract measure is a function $\mu:\mathcal{F}\rightarrow\left[0_\mathbb{S},+\infty\right)$ if the following properties are satisfied:
\begin{enumerate}
\item The abstract measure of any set is a non-negative element of $\mathbb{S}$:
\[
\forall{A\in\mathcal{F}}\;(\mu\left(A\right)\in\mathbb{S})\wedge(\mu\left(A\right)\succeq 0_\mathbb{S});
\]
\item The abstract measure of empty set is equal to $0_\mathbb{S}$: $\mu\left(\emptyset\right)=0_\mathbb{S}$;
\item Monotonicity:  $A\subset B\Rightarrow \mu(A)\preceq\mu(B)$;
\item Countable $\oplus$-additivity: Let $\{A_i\}\subset\mathcal{F}$ is any countable sequence of mutually (pairwise) disjoining sets, then
\[
\mu\left(\bigcup_{i=1}^{\infty}A_i\right)=\bigoplus_{i=1}^{\infty}\mu\left(A_i\right)=\lim_{n\to\infty}\bigoplus_{i=1}^{n}\mu\left(A_i\right).
\]
\end{enumerate}
\end{defi}
An important concept in classical measure theory is the measure space.  For the abstract measure we can also define measure space as measurable space with the abstract measure defined on it. 

Here and below, we always consider $\left(\Omega, \mathcal{F},\mu\right)$ as a measure space with the abstract measure $\mu$.

Let $\left(\Omega, \mathcal{F}\right)$ is a measurable space, then $\Omega$ is considered as the set of all possible outcomes and $\mathcal{F}$ is a set of events that can be composed of the outcomes. So, an event is a specific name of set whose elements are outcomes.

Consider two events $A$ and $B$. Two events can be mutually exclusive, dependent or independent. Moreover, it is obvious that events cannot be mutually exclusive, dependent or independent at the same time.
\begin{defi}
Let $\left(\Omega, \mathcal{F}\right)$ is a measurable space. Exclusivity relation is a homogeneous relation $\parallel$ that is defined over the set $\mathcal{F}$ that $A\parallel B$ if and only if $A\cap B=\emptyset$.
\end{defi}
Then, based on the properties of the intersection of two sets, we can obtain the properties of the exclusivity relation.
\begin{prop}
Let $\left(\Omega, \mathcal{F}\right)$ is a measurable space, then the exclusivity relation $\parallel$ that is defined over the set $\mathcal{F}$ has the following properties:
\begin{enumerate}
\item $\forall A\in \mathcal{F}\neg(A\parallel A)$;
\item $A\parallel B\Leftrightarrow B\parallel A$;
\item $\forall A\in \mathcal{F}(A\parallel\emptyset)$;
\item $(A\parallel C)\wedge(B\parallel C)\Rightarrow (A\cup B)\parallel C$;
\item $(A\parallel C)\wedge(B\parallel C)\Rightarrow (A\cap B)\parallel C$.
\end{enumerate}
\end{prop}
Informally, two events are mutual exclusive if and only if they cannot occur at the same time. In other words, the occurence of one event excludes the occurence of another event. Therefore, we can define mutually exclusive events in terms of exclusivity relation.
\begin{defi}
Two events $A$ and $B$ are mutually exclusive if and only if $A\parallel B$. 
\end{defi}
For countable collection of events we can say that this collection of events is pairwise exclusive if and only if every pair of events are mutually exclusive. Countable collection of events is mutually exclusive if and only if every event is mutually exclusive for every intersection of the other events. Obviously, if collection of events is pairwise exclusive, then it is mutually exclusive, and vice versa.

In the classical aproach to the axiomatization of probability theory presented by Kolmogorov \cite{kol}, the concept of independent events is defined in terms of the probability measure. However, in our opinion, this concept should be defined differently. In fact, we know in advance which events are independent, so it is not the probability measure defines the independent events, but the independent events define the probability measure. 
\begin{defi}
Let $\left(\Omega, \mathcal{F}\right)$ is a measurable space. Independence relation is a homogeneous relation $\perp$ that is defined over the set $\mathcal{F}$ that the following properties are satisfied:
\begin{enumerate}
\item $\forall A\in \mathcal{F}\neg(A\perp A)$;
\item $A\perp B\Leftrightarrow B\perp A$;
\item $\forall A\in \mathcal{F}(A\perp\Omega)$;
\item $(A\perp C)\wedge(B\perp C)\Rightarrow (A\cup B)\perp C$;
\item $(A\perp C)\wedge(B\perp C)\Rightarrow (A\cap B)\perp C$;
\item $A\perp B\Rightarrow A^\mathrm{C}\perp B$.
\end{enumerate}
\end{defi}
Now suppose that $A\cap B\neq \emptyset$, then two events can be mutually independent or dependent. Informally, we can say that two events $A$ and $B$ are independent if and only if the occurrence of $A$ does not affect the abstract probability of occurrence of event $B$, and vice versa. Similarly, two events $A$ and $B$ are dependent if and only if the occurrence of $A$ affect the abstract probability of occurrence of event $B$, and vice versa. Formally, we can define independent and dependent events in terms of the independence relation.
\begin{defi}
\label{ind}
Let $\left(\Omega, \mathcal{F}\right)$ is a measurable space that complex of conditions $\perp$ is defined over the set $\mathcal{F}$. Two events $A$ and $B$ that $A\cap B\neq\emptyset$ are mutually independent if and only if $A\perp B$.
\end{defi}
Consider countable collection of events, then the collection of events is pairwise independent if and only if every pair of events are mutually independent. Similarly, finite collection of events is mutually independent if and only if every event is idependent for every intersection of the other events. 

Note that if collection of events is mutually independent, then the collection is pairwise independent, but the convers is generally not true.

Since events cannot be independent and dependent at the same time, then, based on Definition \ref{ind}, we can define mutually dependent events.
\begin{defi}
\label{dep}
Let $\left(\Omega, \mathcal{F}\right)$ is a measurable space that complex of conditions $\perp$ is defined over the set $\mathcal{F}$. If $A\cap B\neq\emptyset$ and the events are not independent, then they are dependent.
\end{defi}
As we can see, exclusivity and independent relations determine which events are mutually exclusive, independent or dependent.
\begin{defi}
Let $\left(\Omega, \mathcal{F}\right)$ is a measurable space that the exclusivity relation $\parallel$ and the independence relation $\perp$ are defined over the set $\mathcal{F}$, then complex of conditions is a pair $(\parallel,\perp)$.
\end{defi}
Classical probability measure \cite{kol} and possibility measure \cite{dub1}, which we mentioned before, are united by the fact that they take values from the segment $[0,1]$, which can be explained by the algebraic structures that are used. Both algebraic structures have additive identity $0$ and multiplicative identity $1$. However, definitions of these measures do not explain why the probability or the possibility of the sample space is equal to $1$, so in our definition we correct this defect.
\begin{defi}
\label{absprob}
Let $\left(\Omega, \mathcal{F}\right)$ is a measurable space that complex of conditions $(\perp,\parallel)$ is defined over the set $\mathcal{F}$. Abstract probability measure is a function $\mathbb{P}:\mathcal{F}\rightarrow\left[0_\mathbb{S},1_\mathbb{S}\right]$  that the following properties are satisfied:
\begin{enumerate}
\item The abstract probability of any event is non-negative element of $\mathbb{S}$:
\[
\forall{A\in\mathcal{F}}\;(\mathbb{P}\left(A\right)\in\mathbb{S})\wedge(\mathbb{P}\left(A\right)\succeq 0_\mathbb{S});
\]
\item The abstract probability of empty event is equal to $0_\mathbb{S}$: $\mathbb{P}\left(\emptyset\right)=0_\mathbb{S}$;
\item Monotonicity:  $A\subset B\Rightarrow \mathbb{P}(A)\preceq\mathbb{P}(B)$;
\item Countable $\oplus$-additivity: Let $\{A_i\}\subset\mathcal{F}$ is any countable sequence of mutually exclusive events, then
\[
\mathbb{P}\left(\bigcup_{i=1}^{\infty}A_i\right)=\bigoplus_{i=1}^{\infty}\mathbb{P}\left(A_i\right)=\lim_{n\to\infty}\bigoplus_{i=1}^{n}\mathbb{P}\left(A_i\right).
\]
\item The abstract probability of $\Omega$ is equal to $1_\mathbb{S}$: $\mathbb{P}\left(\Omega\right)=1_\mathbb{S}$;
\item Countable $\otimes$-multiplicativity: Let $\{A_i\}\subset\mathcal{F}$ is any countable sequence of mutually independent events, then
\[
\mathbb{P}\left(\bigcap_{i=1}^{\infty}A_i\right)=\bigotimes_{i=1}^{\infty}\mathbb{P}\left(A_i\right)=\lim_{n\to\infty}\bigotimes_{i=1}^{n}\mathbb{P}\left(A_i\right).
\]
\end{enumerate}
\end{defi}
As you can see, the requirment for countable $\otimes$-multiplicativity explains the abstract probability of the sample space is equal to the multiplicative identity. Since for all $A\in\mathcal{F}$ $\Omega\cap A=A$ and, by Definition \ref{ind}, $\Omega\perp A$, then $\mathbb{P}(A\cap\Omega)=\mathbb{P}(A)\otimes\mathbb{P}(\Omega)=\mathbb{P}(A)$. Thus, by Definition \ref{algstr}, $\mathbb{P}(\Omega)=1_\mathbb{S}$.
\begin{defi}
Let $\left(\Omega, \mathcal{F}\right)$ is a measurable space with complex of conditions $(\perp,\parallel)$ and abstract probability measure $\mathbb{P}$ defined on it, then tuple $\left(\Omega, \mathcal{F},\perp,\parallel,\mathbb{P}\right)$ is a probability space.
\end{defi}
Here and below, we always consider $\left(\Omega, \mathcal{F},\perp,\parallel,\mathbb{P}\right)$ as a probability space with the abstract measure $\mathbb{P}$.

It is not difficult to see that the definition of the abstract measure introduced by us is in many ways similar to definitions of probability \cite{kol} and possibility \cite{dub2} measures. Properties of the abstract measure also almost identically repeat the well-known properties that the measures possess.
\begin{theorem}
\label{probprop}
Let $\left(\Omega, \mathcal{F},\perp,\parallel,\mathbb{P}\right)$ is a probability space, then for $A,B\in\mathcal{F}$ the following properties hold:
\begin{enumerate}
\item $A\subseteq B\Rightarrow\mathbb{P}\left(B\setminus A\right)\oplus\mathbb{P}\left(A\right)=\mathbb{P}\left(B\right)$;
\item $\mathbb{P}\left(A\right)\oplus\mathbb{P}\left(A^\mathrm{C}\right)=1_\mathbb{S}$;
\item $1_\mathbb{S}\succeq \mathbb{P}\left(A\right)\succeq 0_\mathbb{S}$;
\item $\mathbb{P}\left(A\cup B\right)\oplus\mathbb{P}\left(A\cap B\right)=\mathbb{P}\left(A\right)\oplus\mathbb{P}\left(B\right)$.
\end{enumerate}
\end{theorem}
\begin{proof}
The properties of the abstract probability measure follow immediately from Definition \ref{algstr} and Definition \ref{absprob}.
\begin{enumerate}
\item Since $A\cup( B\setminus A)=B$, then
\[
\mathbb{P}\left(A\cup( B\setminus A)\right)=\mathbb{P}\left(B\right).
\]
Observe that Definition \ref{absprob} implies that
\[
\mathbb{P}\left(A\cup( B\setminus A)\right)=\mathbb{P}\left(A\right)\oplus\mathbb{P}\left( B\setminus A\right).
\] 
\item Since $A\cup A^\mathrm{C}=\Omega$, then, because of Definition \ref{absprob}, 
\[
\mathbb{P}\left(A\cup A^\mathrm{C}\right)=\mathbb{P}\left(A\right)\oplus\mathbb{P}\left(A^\mathrm{C}\right).
\]
However Definition \ref{absprob} implies that 
\[
\mathbb{P}\left(A\cup A^\mathrm{C}\right)=\mathbb{P}\left(\Omega\right)=1_\mathbb{S}.
\]
Therefore, $\mathbb{P}\left(A\right)\oplus\mathbb{P}\left(A^\mathrm{C}\right)=1_\mathbb{S}$.
\item By Definition \ref{absprob},
\[
\mathbb{P}\left(A\right)\succeq 0_\mathbb{S},\mathbb{P}\left(A^\mathrm{C}\right)\succeq 0_\mathbb{S}.
\]
However note that, by Definition \ref{algstr}, 
\[
\mathbb{P}\left(A^\mathrm{C}\right)\succeq 0\Rightarrow \mathbb{P}\left(A\right)\oplus\mathbb{P}\left(A^\mathrm{C}\right)\succeq \mathbb{P}\left(A\right)\oplus 0_\mathbb{S}.
\]
Therefore, $1_\mathbb{S}\succeq \mathbb{P}\left(A\right)$. That is, $1_\mathbb{S}\succeq \mathbb{P}\left(A\right)\succeq 0_\mathbb{S}$.
\item Firstly, $A\cup B=A\cup\left( B\setminus \left(A\cap B\right)\right)$, therefore, because of Definition \ref{absprob},
\[
\mathbb{P}\left(A\cup B\right)= \mathbb{P}\left(A\cup\left( B\setminus \left(A\cap B\right)\right)\right)= \mathbb{P}\left(A\right)\oplus \mathbb{P}\left( B\setminus \left(A\cap B\right)\right).
\] 
Secondly, by Definition \ref{algstr}, 
\[
\mathbb{P}\left(A\cup B\right)\oplus\mathbb{P}\left(A\cap B\right)=\mathbb{P}\left(A\right)\oplus \mathbb{P}\left( B\setminus \left(A\cap B\right)\right)\oplus\mathbb{P}\left(A\cap B\right).
\]
Thus, $\mathbb{P}\left(A\cup B\right)\oplus\mathbb{P}\left(A\cap B\right)=\mathbb{P}\left(A\right)\oplus\mathbb{P}\left(B\right)$.
\end{enumerate}
\end{proof}
In addition to all of the above, the abstract probability measure can be used to describe the abstract probability distribution of random variable.
\begin{defi}
\label{probdist}
Let $\left(\Omega,\mathcal{F},\parallel_\mathcal{F},\perp_\mathcal{F},\mathbb{P}\right)$ and $\left(E, \mathcal{E},\parallel_\mathcal{E},\perp_\mathcal{E},X_{*}\mathbb{P}\right)$ are  probability spaces and $X:\Omega\to E$ is a $\left(E,\mathcal{E}\right)$-valued random variable, then an abstract probability distribution of the random variable is the abstract probability measure $ X_{*}\mathbb{P}:\mathcal{E}\rightarrow\left[0_\mathbb{S},1_\mathbb{S}\right]$ pushforward of the abstract probability measure $\mathbb{P}:\mathcal{F}\rightarrow\left[0_\mathbb{S},1_\mathbb{S}\right]$ that
\[
\forall A\in\mathcal{E} \left( X_{*}\mathbb{P}\left(A\right)=\mathbb{P}\left(X^{-1}(A)\right)\right).
\]
\end{defi}
Since abstract probability distribution is abstract probability measure defined on measurable space $(E,\mathcal{E})$, then it satisfies axioms presented in Definition \ref{absprob}. Therefore, the definition implies that abstract probability distribution  has the properties presented in Theorem \ref{probprop}.

Based on Definition \ref{ind} and Definition \ref{probdist}, we can define mutually independent variables.
\begin{defi}
Let $\left(\Omega,\mathcal{F},\parallel_\mathcal{F},\perp_\mathcal{F},\mathbb{P}\right)$ and $\left(E, \mathcal{E},\parallel_\mathcal{E},\perp_\mathcal{E},X_{*}\mathbb{P}\right)$ are  probability spaces and $X,Y:\Omega\to E$ are $\left(E,\mathcal{E}\right)$-valued random variables. Two variables $X$ and $Y$ are mutually independent if and only if events $X^{-1}(A)$ and $Y^{-1}(B)$  are independent for all $A,B\in\mathcal{E}$.
\end{defi}
By similar reasoning, we can define pairwise and mutual independent random variables. Let $\{X_i\}_{i\in I}$ is a finite collection of $(E,\mathcal{E})$-valued random variables defined on  probability space $\left(\Omega,\mathcal{F},\parallel,\perp,\mathbb{P}\right)$, then the collection of random variables is pairwise independent if and only if for all collections of events $\{A_i\}_{i\in I}\subset\mathcal{E}$ collection of events $\{X_i^{-1} (A_i)\}_{i\in I}$ is pairwise independent, where $I$ is an index set. The collection of random variables is mutually independent if and only if for all collections of events $\{A_i\}_{i\in I}\subset\mathcal{E}$ collection of events $\{X_i^{-1}( A_i)\}_{i\in I}$ is mutually independent, where $I$ is an index set.
\section{Abstract conditional probability}
Informally, we can say that conditional probability as a specific abstract probability measure of event $A$, given the occurrence of event $B$.

Consider two events $A$ and $B$. If two events are mutually exclusive, then abstract conditional probability of event $A$ given the occurence of event $B$  is equal to the abstract probability of  impossible event. If $A$ and $B$ are independent, then abstract conditional probability of event $A$ given the occurence of event $B$ is an abstract probability of event $A$. Similarly, if $A$ and $B$ are independent, then abstract conditional probability of event A given the occurence of event B is not equal to an abstract probability of event A.
\begin{defi}
\label{condprob}
Let $\left(\Omega,\mathcal{F},\parallel,\perp,\mathbb{P}\right)$ and $\left(\Omega,\mathcal{F},\parallel,\perp,\mathbb{P}_B\right)$ are  probability spaces, $A,B\in\mathcal{F}$ are events that $\mathbb{P}(B)\neq 0_\mathbb{S}$. Abstract conditional probability is an abstract probability measure $\mathbb{P}_B:\mathcal{F}\to[0_\mathbb{S},1_\mathbb{S}]$ that the following properties are satisfied:
\begin{enumerate}
\item $A\parallel B\Rightarrow\mathbb{P}_B(A)=0_\mathbb{S}$;
\item $A\perp B\Rightarrow\mathbb{P}_B(A)=\mathbb{P}(A)$;
\item $\neg(A\perp B)\Rightarrow\mathbb{P}_B(A)\neq\mathbb{P}(A)$;
\item $\mathbb{P}\left(A\cap B\right)=\mathbb{P}_B\left(A\right)\otimes\mathbb{P}\left(B\right)$.
\end{enumerate}
\end{defi}
In addition, from Definition \ref{condprob} immediately follows that abstract conditional probability can be used as a criterion for defining independent events.
\begin{theorem}
Two events $A$ and $B$ are mutually independent if and only if $\mathbb{P}_B(A)=\mathbb{P}(A)$.
\end{theorem}
\begin{proof}
Firstly, by Definition \ref{dep}, two events $A$ and $B$ are mutually dependent if and only if $A\cap B\neq \emptyset$ and $\neg(A\perp B)$.

Secondly, by Definition \ref{condprob}, if $A\cap B\neq \emptyset$, then $(A\perp B)\rightarrow(\mathbb{P}_B(A)=\mathbb{P}(A))$ and $\neg(A\perp B)\rightarrow\neg(\mathbb{P}_B(A)=\mathbb{P}(A))$.

Therefore, since $\forall a,b (a\leftrightarrow b=(a\rightarrow b)\wedge(\neg a\rightarrow\neg b))$, then
\[
((A\perp B)\rightarrow(\mathbb{P}_B(A)=\mathbb{P}(A)))\wedge(\neg(A\perp B)\rightarrow\neg(\mathbb{P}_B(A)=\mathbb{P}(A)))=(A\perp B)\leftrightarrow(\mathbb{P}_B(A)=\mathbb{P}(A)).
\]

Thus, $A\perp B\Leftrightarrow\mathbb{P}_B(A)=\mathbb{P}(A)$.
\end{proof}
Since the abstract conditional probability is the abstract probability measure, then it has the properties presented in Definition \ref{absprob} and Theorem \ref{probprop}. The abstract conditional probability also has analogues of the other well-recognized properties of classical conditional probability as the Bayes' theorem and the law of total probability described below.
\begin{theorem}
Let $\left(\Omega,\mathcal{F},\parallel,\perp,\mathbb{P}\right)$ and $\left(\Omega,\mathcal{F},\parallel,\perp,\mathbb{P}_B\right)$ are  probability spaces, $A,B\in\mathcal{F}$ are events that $\mathbb{P}(A)\neq 0_\mathbb{S}$ and $\mathbb{P}(B)\neq 0_\mathbb{S}$, then  
\[
\mathbb{P}\left(A\right)\otimes\mathbb{P}_A\left( B\right)=\mathbb{P}_B\left(A\right)\otimes\mathbb{P}\left( B\right).
\]
\end{theorem}
\begin{proof}
Observe that, if $(\mathbb{S},\oplus,\otimes,0_\mathbb{S},1_\mathbb{S},\succeq)$ is an ordered pseudo-semiring with multiplicative inverse, then Definition \ref{condprob} implies that 
\[
\mathbb{P}\left(A\cap B\right)=\mathbb{P}\left(A\right)\otimes\mathbb{P}_A\left( B\right)
\]
and
\[
\mathbb{P}\left(B\cap A\right)=\mathbb{P}_B\left(A\right)\otimes\mathbb{P}\left( B\right).
\]
Therefore, since $A\cap B=B\cap A$ 
\[
\mathbb{P}\left(A\right)\otimes\mathbb{P}_A\left( B\right)=\mathbb{P}_B\left(A\right)\otimes\mathbb{P}\left( B\right).
\]
\end{proof}
\begin{theorem}
Let $\left(\Omega,\mathcal{F},\parallel,\perp,\mathbb{P}\right)$ and $\left(\Omega,\mathcal{F},\parallel,\perp,\mathbb{P}_B\right)$ are  probability spaces, $A,B\in\mathcal{F}$ are events that $\mathbb{P}(A)\neq 0_\mathbb{S}$ and $\mathbb{P}(B)\neq 0_\mathbb{S}$. If $A\in\mathcal{F}$ is any event and $\{H_i\}_{i\in I}\mathcal{F}$ is a finite set of pairwise exclusive events that $\bigcup_{i\in I}H_i=\Omega$ and $\forall i\in I(\mathbb{P}(H_i)\neq 0_\mathbb{S})$, then
\[
\mathbb{P}\left(A\right)=\bigoplus_{i\in I}\mathbb{P}\left(H_i\right)\otimes\mathbb{P}_{H_i}\left(A\right),
\]
where $I$ is an index set.
\end{theorem}
\begin{proof}
Firstly, if $\{H_i\}_{i\in I}\mathcal{F}$ is a finite set of pairwise exclusive events that $\bigcup_{i\in I}H_i=\Omega$, then
\[
A=A\cap\Omega=A\cap\left(\bigcup_{i\in I}H_i\right)=\bigcup_{i\in I}\left(A\cap H_i\right).
\]
Since $\{H_i\}_{i\in I}$ is a finite set of pairwise exclusive events, then $\{A\cap H_i\}_{i\in I}$ is also a finite set of pairwise exclusive events. Therefore, because of Definition \ref{absprob}, 
\[
\mathbb{P}\left(A\right)=\mathbb{P}\left(\bigcup_{i\in I}\left(A\cap H_i\right)\right)=\bigoplus_{i\in I}\mathbb{P}\left(A\cap H_i\right).
\]
Thus, by Definitiom \ref{condprob}, if $\forall i\in I(\mathbb{P}(H_i)\neq 0_\mathbb{S})$, then
\[
\mathbb{P}\left(A\right)=\bigoplus_{i\in I}\mathbb{P}\left(H_i\right)\otimes\mathbb{P}_{H_i}\left(A\right).
\]
\end{proof}
\section{Abstract expected value operator as abstract Lebesgue integral}
Classical measure has found wide application in the Lebesgue integration of functions \cite{burkill}, where it plays key role. 

As generalization of the Lebesgue integral, authors have presented various definitions of integral that are in many ways similar to the Lebesgue integral, but instead of the classical measure, they use the corresponding generalized measures that we mentioned before. Examples of integrals are Sugeno's fuzzy integral \cite{sugeno}, pseudo-integral \cite{sugeno1}, generalized Lebesgue integral \cite{zhang}, etc. 

All the generalized integrals described above are united by fact that they can integrate only real-valued functions. In this section, we define an abstract Lebesgue integral that can be used to integrate a wider class of functions.
\begin{defi}
\label{simplef}
Let $(\Omega,\mathcal{F})$ is a measurable space, $\{a_i\}_{i=1}^n\subset\mathbb{S}$ is a sequence of elements of $\mathbb{S}$ and $\{A_i\}_{i=1}^n\subset\mathcal{F}$ is a sequence of a pairwise disjoint measurable sets that $\bigcup_{i=1}^{n}A_i=\Omega$, then abstract simple function is a function $f:\Omega\to\mathbb{S}$ such that
\[
f(\omega)=\bigoplus_{i=1}^{n}a_i\otimes \chi_{A_i}(\omega)
\]
where 
\[
\chi_{A_i}(\omega)=\begin{cases}1_\mathbb{S}, &\omega\in A_i\\0_\mathbb{S},&\omega\notin A_i.\end{cases}
\]
\end{defi}
In other words, if $f(\omega)=\bigoplus_{i=1}^{n}a_i\otimes \chi_{A_i}(\omega)$, then since $0_\mathbb{S}$ is additive identity and $1_\mathbb{S}$ is multiplicative identity $\forall \omega\in A_i (f(\omega)=a_i)$.
\begin{defi}
\label{integral}
Let $\left(\Omega, \mathcal{F},\mu\right)$ is a measure space. 
\begin{enumerate}
\item If $f:\Omega\to\mathbb{S}$ is an abstract simple function that
\[
f(\omega)=\bigoplus_{i=1}^{n}a_i\otimes \chi_{A_i}(\omega),
\]
then an abstract Lebesgue integral of the abstract simple function is given by
\[
\bigoplus_A f(\omega)\otimes\mathrm{d}\mu(\omega)=\bigoplus_{i=1}^{n}a_i\otimes \mu(A_i\cap A).
\]
\item Let $\{f_n\}$ is a sequence of abstract simple functions that $f_n:\Omega\to\mathbb{S}$ and $f_n\nearrow f$ as $n\to\infty$, then an abstract Lebesgue integral of the measurable function $f:\Omega\to\mathbb{S}$ is given by
\[
\bigoplus_A f(\omega)\otimes\mathrm{d}\mu(\omega)=\lim_{n\to\infty}\bigoplus_A f_n(\omega)\otimes\mathrm{d}\mu(\omega).
\]
\end{enumerate}
\end{defi}
Let  $\left(\Omega, \mathcal{F},\mu\right)$ is a measure space, then any measurable function $f:\Omega\to\mathbb{S}$ can be represented as follows:
\[
f(\omega)=f^\oplus(\omega)\oplus f^\ominus(\omega)
\]
where $f^\oplus(\omega)=\max\{0_\mathbb{S},f(\omega)\}$ is a non-negative part and $f^\ominus(\omega)=\min\{0_\mathbb{S},f(\omega)\}$ is a non-positive part. 

We say that the abstract Lebesgue integral of the measurable function exists if at least one of integrals $\bigoplus_A f^\oplus(\omega)\otimes\mathrm{d}\mu(\omega)$ and $\bigoplus_A f^\ominus(\omega)\otimes\mathrm{d}\mu(\omega)$ is finite.

As you can see, the definition of the abstract Lebesgue integral introduced by us is in many ways similar to definitions of Lebesgue integral \cite{burkill}, the generalized Lebesgue integral \cite{zhang} and pseudo-integral \cite{zhang1}. Properties of the abstract Lebesgue integral also almost identically repeat the well-known properties that the mentioned integrals possess.
\begin{theorem}
\label{intprop}
Let $\left(\Omega, \mathcal{F},\mu\right)$ is a measure space, then for any $A\in\mathcal{F}$ an abstract Lebesgue integral has the following properties:
\begin{enumerate}
\item If $f,g:\Omega\to\mathbb{S}$ are measurable functions that $f(\omega)\preceq g(\omega)$, then 
\[
\bigoplus_{A}f(\omega)\otimes\mathrm{d}\mu(\omega)\preceq\bigoplus_{A}g(\omega)\otimes\mathrm{d}\mu(\omega);
\]
\item Let $f:\Omega\to\mathbb{S}$ is a measurable function, then for all $a\in\mathbb{S}$
\[
\bigoplus_{A}a\otimes\mathrm{d}\mu(\omega)=a\otimes\mu(A);
\]
\item If $f,g:\Omega\to\mathbb{S}$ are measurable functions, then
\[
\bigoplus_{A}(f(\omega)\oplus g(\omega))\otimes\mathrm{d}\mu(\omega)=\bigoplus_{A}f(\omega)\otimes\mathrm{d}\mu(\omega)\oplus\bigoplus_{A}g(\omega)\otimes\mathrm{d}\mu(\omega);
\]
\item Let $f:\Omega\to\mathbb{S}$ is a measurable function, then for all $a\in\mathbb{S}$
\[
\bigoplus_{A}a\otimes f(\omega)\otimes\mathrm{d}\mu(\omega)=a\otimes\bigoplus_{A}f(\omega)\otimes\mathrm{d}\mu(\omega);
\]
\end{enumerate}
\end{theorem}
\begin{proof}
The properties of the abstract Lebesgue integral follow immediately from Definition \ref{simplef} and Definition \ref{integral}.
\begin{enumerate}
\item Let $f,g:\Omega\to\mathbb{S}$ are measurable functions that $f\preceq g$, then, by Definition \ref{algstr}, there are sequences of abstract simple functions $\{f_n\}$ and $\{g_n\}$, where  $f_n,g_n:\Omega\to\mathbb{S}$, that $\forall n(f_n\preceq g_n)$, $f_n\nearrow f$ and $g_n\nearrow g$ as $n\to\infty$. Observe that, since $f_n\preceq g_n$, then 
\[
\bigoplus_A f_n(\omega)\otimes\mathrm{d}\mu(\omega)\preceq\bigoplus_A g_n(\omega)\otimes\mathrm{d}\mu(\omega).
\]
Therefore, by Definition \ref{algstr},
\[
\lim_{n\to\infty}\bigoplus_A f_n(\omega)\otimes\mathrm{d}\mu(\omega)\preceq\lim_{n\to\infty}\bigoplus_A g_n(\omega)\otimes\mathrm{d}\mu(\omega).
\]
Thus, by Definition \ref{integral},
\[
\bigoplus_{A}f(\omega)\otimes\mathrm{d}\mu(\omega)\preceq\bigoplus_{A}g(\omega)\otimes\mathrm{d}\mu(\omega).
\]
\item If $f:\Omega\to\mathbb{S}$ is an abstract simple function that
\[
f(\omega)=\bigoplus_{i=1}^{n}a\otimes \chi_{A_i}(\omega)=a,
\]
then, by Definition \ref{integral}, an abstract Lebesgue integral of the abstract simple function is given by
\[
\bigoplus_A f(\omega)\otimes\mathrm{d}\mu(\omega)=\bigoplus_{i=1}^{n}a\otimes \mu(A_i\cap A).
\]
However, note that, because of Definition \ref{algstr} and Definition \ref{meas},
\[
\bigoplus_{i=1}^{n}a\otimes \mu(A_i\cap A)=a\otimes\bigoplus_{i=1}^{n} \mu(A_i\cap A)=a\otimes\mu(A).
\]
Thus, 
\[
\bigoplus_A a\otimes\mathrm{d}\mu(\omega)=a\otimes\mu(A).
\]
\item Let $f,g:\Omega\to\mathbb{S}$ are measurable functions, $\{f_n\}$ and $\{g_n\}$ are sequences of abstract simple functions, where  $f_n,g_n:\Omega\to\mathbb{S}$, that $f_n\nearrow f$ and $g_n\nearrow g$ as $n\to\infty$, then, by Definition \ref{algstr}, $f_n\oplus g_n\nearrow f\oplus g$ as $n\to\infty$. Therefore, by Definition \ref{algstr} and Definition \ref{integral},
\[
\begin{split}
\bigoplus_{A}(f(\omega)\oplus g(\omega))\otimes\mathrm{d}\mu(\omega)&=\lim_{n\to\infty}\bigoplus_A f_n(\omega)\oplus g_n(\omega)\otimes\mathrm{d}\mu(\omega)\\&=\lim_{n\to\infty}\left(\bigoplus_A f_n(\omega)\otimes\mathrm{d}\mu(\omega)\oplus\bigoplus_A g_n(\omega)\otimes\mathrm{d}\mu(\omega)\right)\\&=\lim_{n\to\infty}\bigoplus_A f_n(\omega)\otimes\mathrm{d}\mu(\omega)\oplus\lim_{n\to\infty}\bigoplus_A g_n(\omega)\otimes\mathrm{d}\mu(\omega)\\&=\bigoplus_{A}f(\omega)\otimes\mathrm{d}\mu(\omega)\oplus\bigoplus_{A}g(\omega)\otimes\mathrm{d}\mu(\omega).
\end{split}
\]
\item Let $f:\Omega\to\mathbb{S}$ is measurable function, $\{f_n\}$ is sequence of abstract simple functions, where  $f_n:\Omega\to\mathbb{S}$, that $f_n\nearrow f$ as $n\to\infty$, then, by Definition \ref{algstr}, $\forall a\in\mathbb{S}(a\otimes f_n\nearrow a\otimes f)$ as $n\to\infty$. Therefore, by Definition \ref{algstr} and Definition \ref{integral},
\[
\begin{split}
\bigoplus_{A}a\otimes f(\omega)\otimes\mathrm{d}\mu(\omega)&=\lim_{n\to\infty}\bigoplus_A a\otimes f_n(\omega)\otimes\mathrm{d}\mu(\omega)\\&=\lim_{n\to\infty}\left(a\otimes\bigoplus_A f_n(\omega)\otimes\mathrm{d}\mu(\omega)\right)\\&=a\otimes\lim_{n\to\infty}\bigoplus_A f_n(\omega)\otimes\mathrm{d}\mu(\omega)\\&=a\otimes\bigoplus_{A}f(\omega)\otimes\mathrm{d}\mu(\omega).
\end{split}
\]
\end{enumerate}
\end{proof}
In the classical aproach to the axiomatization of probability theory presented by Kolmogorov \cite{kol}, the expected value of a real-valued random variable is defined as the Lebesgue integral of the random variable as measurable function. Therefore, by similar reasoning, we can define an abstract expected value of $\left(\mathbb{S},\mathcal{S}\right)$-valued random variable in terms of the abstract Lebesgue integral.
\begin{defi}
Let $X$ is a measurable $\left(\mathbb{S},\mathcal{S}\right)$-valued random variable defined on a probability space $\left(\Omega, \mathcal{F},\perp,\parallel,\mathbb{P}\right)$, then an abstract expected value operator of the random variable is defined as the abstract Lebesgue integral of a measurable function
\[
\mathbb{E}[X]=\bigoplus_{\Omega}X(\omega)\otimes \mathrm{d}\mathbb{P}(\omega).
\]
\end{defi}
Since the abstract expected value operator is defined as the abstract Lebesgue integral, then the properties of the operator follow immediately from Theorem \ref{intprop}.
\begin{prop}
Let $X$ and $Y$ are measurable $\left(\mathbb{S},\mathcal{S}\right)$-valued random variables defined on a probability space $\left(\Omega, \mathcal{F},\perp,\parallel,\mathbb{P}\right)$, then the abstract expected value operator has the following properties:
\begin{enumerate}
\item Monotonicity: $X\preceq Y\Rightarrow\mathbb{E}[X]\preceq\mathbb{E}[Y]$;
\item Expectation of constant: $\forall \omega\in \Omega(X=a)\Rightarrow\mathbb{E}[a]=a$;
\item Linearity: $\forall a,b\in\mathbb{S}\left(\mathbb{E}[a\otimes X\oplus b\otimes Y]=a\otimes\mathbb{E}[X]\oplus b\otimes\mathbb{E}[Y]\right)$.
\end{enumerate}
\end{prop}
\section{Conclusions}
We have presented all the main components of the abstract probability theory: abstract probability measure, abstract conditional probability measure and abstract expected value operator.

For the purpose of axiomatization of the abstract probability theory, were defined abstract measure and abstract Lebesgue integral on which our study is mainly based on. 

In addition, we revised and criticized the Kolmogorov's system of axioms \cite{kol}. In order to correct the defect was presented the complex of conditions, which was later used in the definition of the abstract probability measure. 

We have also proven the properties of the abstract probability measure, abstract conditional measure and abstract expected value. The properties we obtained are abstract analogues of the corresponding properties presented in the classical probability theory, in particular , the sum rule, product rule, Bayes' theorem, law of total probability, etc.

Moreover, it is not hard to see that the Kolmogorov's probability theory \cite{kol} and possibility theory \cite{dub1} are special cases of the abstract probability theory.

As a result, we can say that the theory presented by us is not only inferior in anything to the generally accepted Kolmogorov's probability theory, but also in many ways surpasses it, since it has wider application.

\end{document}